\documentclass[11pt]{article}
\usepackage[final]{epsfig}
\usepackage[left=3cm,right=3cm, top=2.5cm,bottom=2.5cm,bindingoffset=0cm]{geometry}
\usepackage{graphics}
\usepackage{amsmath}
\usepackage{amsfonts}
\usepackage{latexsym}
\usepackage{amssymb, mathabx}
\usepackage{amsthm}
\usepackage{graphicx}
\usepackage{epstopdf}
\usepackage{multicol,multirow}
\usepackage{hyperref, enumitem}
\usepackage{mathtools}

\usepackage{mathrsfs}

 \usepackage{relsize}

\usepackage[nottoc]{tocbibind}

\usepackage{amsthm}

\usepackage[latin1]{inputenc}
\usepackage{tikz}
\usetikzlibrary{shapes,arrows,decorations.pathmorphing}
\usepackage{float}

\usepackage{blkarray}
\usepackage{xypic}

\usepackage[title]{appendix}

\numberwithin{equation}{section}

\makeatletter
\def\thmhead@plain#1#2#3{%
  \thmname{#1}\thmnumber{\@ifnotempty{#1}{ }\@upn{#2}}%
  \thmnote{ {\the\thm@notefont#3}}}
\let\thmhead\thmhead@plain

\makeatother

\newcounter{AppCounter}

\def\restrict#1{\raise-.5ex\hbox{\ensuremath|}_{#1}}

\newtheorem{lemma}{Lemma}[section]
\newtheorem{proposition}[lemma]{Proposition}

\newtheorem{remark-definition}[lemma]{Remark-Definition}
\newtheorem{theorem}[lemma]{Theorem}

\newtheorem{proposition-conjecture}[lemma]{Proposition-conjecture}

\theoremstyle{definition}

\newtheorem{remark}[lemma]{Remark}

\newcommand{\marginnote}[1]
{
}

\newcounter{cy}

\newcounter{bk}

\newcounter{dps}

\usepackage{lipsum}

\title{On the Finite-Time Singularities of Point-Vortex Systems II: Upper Right Quadrant}

\author{ Cheng Yang\thanks{Division of Mathematical Sciences, Nanyang Technological University, 21 Nanyang Link, Singapore 637371; 
 e-mail: \tt{cheng.yang@ntu.edu.sg} }
\\
}

\begin{document}
\date{}
\maketitle

\begin{abstract}
We study finite-time singularities of point-vortex systems
in the upper right quadrant. For two vortices with arbitrary
nonzero circulations and distinct initial positions, we
prove global existence and uniqueness by excluding interior
collisions, boundary approach, and finite-time escape to
infinity. For three vortices with arbitrary nonzero
circulations, we rule out finite-time boundary approach
in the limit inferior sense, without assuming boundedness
or convergence of the trajectories. The argument combines
Hamiltonian estimates, monotonicity of signed coordinate
moments, and geometric estimates adapted to both walls
and the corner.

We also establish the existence of finite-time interior
collapse for every nonzero circulation triple whose pairwise
products sum to zero. The collapse point can be prescribed
arbitrarily in the interior, with any sufficiently small
positive collapse time. The resulting solutions are
asymptotically self-similar, exhibiting square-root
contraction and logarithmic rotation. They are constructed
using a local strong stable manifold of real dimension
two in a similarity system augmented by a complex scale.
\end{abstract}

\tableofcontents


\section{Introduction}
\label{sec:introduction}

Point-vortex dynamics provides a finite-dimensional model of
two-dimensional incompressible ideal flow in which vorticity
is concentrated at finitely many moving points. Although
the governing equations have a classical Hamiltonian
structure, their singular interactions permit complicated
dynamics and finite-time collisions. Understanding when
such singularities occur, and how the geometry of the fluid
domain affects them, remains a central problem in vortex
dynamics; see, for example, \cite{Newton2001}.

In the full plane, translation and rotation invariance yield
conservation of linear and angular impulse. These quantities,
together with the Hamiltonian, play an important role in the
analysis of the three-vortex problem. A classical example of
singularity formation is the self-similar collapse discovered
by Gr\"obli: for suitable mixed-sign circulations and initial
configurations, a vortex triangle contracts and rotates until
its three vertices meet in finite time. The circulation
condition
\begin{equation}
\Gamma_1\Gamma_2+\Gamma_1\Gamma_3+\Gamma_2\Gamma_3=0
\label{eq:intro-collapse-condition}
\end{equation}
allows such collapsing configurations, whose characteristic
length scale is proportional to the square root of the
remaining time. The geometry and stability of these motions
are discussed in \cite{Aref2010}.

For domains with fewer continuous symmetries than the full
plane, including the half-plane and the quadrant, the
classification of finite-time singularities has long
presented difficult open questions
\cite{Newton2001,Yang}. Impermeable boundaries introduce
additional interactions through the Green function and its
regular part, as described by the Kirchhoff--Routh
Hamiltonian \cite{Lin1941}. These boundary-induced terms
become singular near the boundary and remove some of the
conservation laws available in the plane. In particular,
energy conservation alone may not prevent boundary approach: 
the energy from vortex interactions can increase without bound 
while the energy from boundary interactions decreases, keeping the total energy constant.

Several results clarify different aspects of this problem.
Donati \cite{Donati2022} proved global existence for almost
every initial configuration in bounded planar domains with
\(C^{2,\alpha}\) boundary, leaving open the behavior of
exceptional initial data. Donati, Godard-Cadillac, and Iftimie
\cite{DoGoIf} excluded finite-time
boundary collapse for positive circulations in the half-plane
and disk.  In Part~I of this series
\cite{Yang2}, we excluded finite-time boundary
approach for three vortices with arbitrary nonzero
circulations in the half-plane and constructed interior
collapsing solutions.

The present paper considers the upper right quadrant
\[
\Omega=\{z=x+iy\in\mathbb C:x>0,\ y>0\}.
\]
This geometry introduces two interacting walls and their
common corner. Although the method of images gives an
explicit velocity field, the tangential translation symmetry
of the half-plane is lost, and the corner requires estimates
that remain valid near both walls simultaneously. Moreover,
the domain is unbounded, so boundary exclusion must account
for possible tangential escape and repeated excursions,
rather than assuming convergence of the trajectories.

\subsection*{Main Results and Organization of the Paper}

Our main results concern the exclusion and formation of
finite-time singularities for two and three vortices in
the quadrant.

For two vortices with arbitrary nonzero circulations, we
first exclude simultaneous accumulation at every finite
point of \(\overline{\Omega}\), whether in the interior,
on either wall, or at the corner (Theorems~\ref{thm:two_interior} and \ref{thm:no-boundary-accumulation}). More precisely, for every
classical solution with distinct positions on
\([0,T)\), \(T<\infty\), we prove
\[
\liminf_{t\uparrow T}
\bigl(|z_1(t)-p|+|z_2(t)-p|\bigr)>0
\qquad\text{for every }p\in\overline{\Omega}.
\]
Together with estimates excluding finite-time escape to
infinity and solitary boundary approach, these results
establish global existence and uniqueness for every pair
of distinct initial positions and every choice of nonzero
circulations (Theorem~\ref{thm:two_global}). The solution exists for all
\(t\in\mathbb R\). On every bounded time interval, the
positions remain bounded, and both the mutual distance
and the distances to the boundary have positive lower bounds.

For three vortices, we prove finite-time boundary exclusion
without any restriction on the circulation signs.
Theorem~\ref{thm:three-boundary-exclusion-unbounded} states
that every classical solution with pairwise distinct
positions on \([0,T)\), \(T<\infty\), satisfies
\[
\inf_{0\le t<T}\min_{1\le j\le3}
\operatorname{dist}(z_j(t),\partial\Omega)>0.
\]
No boundedness or trajectory-convergence assumption is
required. In particular, neither wall nor the corner can
be approached in the limit inferior sense, even when
unbounded excursions or nonconvergent motions are allowed.

Interior collapse nevertheless occurs for suitable
mixed-sign initial configurations.
Theorem~\ref{thm:collapse} applies to every nonzero
circulation triple satisfying
\[
\Gamma_1\Gamma_2+\Gamma_1\Gamma_3+\Gamma_2\Gamma_3=0.
\]
For every prescribed \(z_*\in\Omega\) and every sufficiently
small collapse time \(t_*>0\), we construct a solution
whose three vortices remain distinct before \(t_*\) and
satisfy
\[
z_j(t)
=
z_*+(t_*-t)^{1/2+i\mu}A_j+O(t_*-t),
\qquad j=1,2,3.
\]
Here \(\mathbf A\) is a noncollinear profile with pairwise
distinct components, and \(\mu\ne0\); both depend only on
the circulations. The trajectories remain in a fixed
interior neighborhood of \(z_*\). These solutions are
asymptotically self-similar, with square-root contraction
and logarithmic rotation. 
Their construction uses a local strong stable manifold
of real dimension two for the similarity system augmented
by a complex scale, which encodes both contraction and
logarithmic rotation.

The remainder of the paper is organized as follows.
Section~\ref{sec: pre} introduces the equations for \(N\)
point vortices in the upper right quadrant, together with
their Hamiltonian and complex formulations.
Section~\ref{sec:Nvor} establishes basic estimates for the
\(N\)-vortex system. In particular, it shows that a vortex
approaching either wall or the corner in finite time must
itself come arbitrarily close to another vortex. It also
excludes interior collapse when \(N=2\).
Section~\ref{sec:BoundaryNoncollisionTwo} rules out
simultaneous boundary accumulation for two vortices.
Combining this result with interior noncollision and
estimates excluding finite-time escape to infinity, we
establish global existence for two vortices with arbitrary
nonzero circulations and distinct initial positions.
Section~\ref{sec:three-boundary-exclusion} excludes
finite-time boundary approach for three vortices, without
assuming bounded or convergent trajectories.
Finally, Section~\ref{sec:rigorous_existence_collapse}
constructs finite-time interior collapsing solutions
in the quadrant.

\medskip


\section{Preliminaries}\label{sec: pre}

Consider the first quadrant $\Omega = \mathbb{R}^2_{++}= \{(x,y) \in \mathbb{R}^2 \mid x > 0, y > 0\}$. For $N$ point vortices with positions $(x_j, y_j) \in \Omega$ and circulations $\Gamma_j$, impermeability along $x=0$ and $y=0$ is enforced via the method of images. Each vortex at $z_k = x_k + i y_k$ induces three image vortices: $(x_k, -y_k)$ with circulation $-\Gamma_k$ (reflection across $x$-axis), $(-x_k, y_k)$ with circulation $-\Gamma_k$ (reflection across $y$-axis), and $(-x_k, -y_k)$ with circulation $+\Gamma_k$ (reflection through origin).

The motion is governed by the Hamiltonian $H$ via $\Gamma_j \dot{x}_j = \frac{\partial H}{\partial y_j}$ and $\Gamma_j \dot{y}_j = -\frac{\partial H}{\partial x_j}$:
\begin{equation}\label{eq:1/4_Ham}
H = \sum_{1 \le j < k \le N} \Gamma_j \Gamma_k G(x_j, y_j; x_k, y_k) + \frac{1}{8\pi} \sum_{j=1}^N \Gamma_j^2 \ln \left( \frac{4 x_j^2 y_j^2}{x_j^2 + y_j^2} \right),
\end{equation}
where $G(x_j, y_j; x_k, y_k)$ is the Dirichlet Green's function for the quadrant:
\begin{equation}\label{eq:1/4_Green}
G = -\frac{1}{4\pi} \ln \left[ \frac{\big((x_j-x_k)^2 + (y_j-y_k)^2\big)\big((x_j+x_k)^2 + (y_j+y_k)^2\big)}{\big((x_j-x_k)^2 + (y_j+y_k)^2\big)\big((x_j+x_k)^2 + (y_j-y_k)^2\big)} \right].
\end{equation}
For each vortex $j = 1, \dots, N$, define the four squared distances to vortex $k$ and its images:
$$d_{1, jk}^2 = (x_j - x_k)^2 + (y_j - y_k)^2, \quad d_{2, jk}^2 = (x_j - x_k)^2 + (y_j + y_k)^2,$$
$$d_{3, jk}^2 = (x_j + x_k)^2 + (y_j - y_k)^2, \quad d_{4, jk}^2 = (x_j + x_k)^2 + (y_j + y_k)^2.$$
The equations of motion split into self-interaction (boundary bounce) and inter-vortex terms:
\begin{equation}\label{eq:1/4_vortex_eq}
\begin{aligned}
\dot{x}_j &= \frac{\Gamma_j}{4\pi} \frac{x_j^2}{y_j (x_j^2 + y_j^2)} + \frac{1}{2\pi} \sum_{\substack{k=1 \\ k \neq j}}^N \Gamma_k \left[ -\frac{y_j - y_k}{d_{1, jk}^2} + \frac{y_j + y_k}{d_{2, jk}^2} + \frac{y_j - y_k}{d_{3, jk}^2} - \frac{y_j + y_k}{d_{4, jk}^2} \right],\\
\dot{y}_j &= -\frac{\Gamma_j}{4\pi} \frac{y_j^2}{x_j (x_j^2 + y_j^2)} + \frac{1}{2\pi} \sum_{\substack{k=1 \\ k \neq j}}^N \Gamma_k \left[ \frac{x_j - x_k}{d_{1, jk}^2} - \frac{x_j - x_k}{d_{2, jk}^2} - \frac{x_j + x_k}{d_{3, jk}^2} + \frac{x_j + x_k}{d_{4, jk}^2} \right].
\end{aligned}
\end{equation}

Using complex velocity $\frac{d\bar{z}_j}{dt} = \dot{x}_j - i \dot{y}_j$, the equations can also be written as:
\begin{equation}\label{eq:1/4_vortex_complex}
\frac{d\bar{z}_j}{dt} = \frac{\Gamma_j}{2\pi i} \left( \frac{1}{2 z_j} - \frac{1}{z_j - \bar{z}_j} - \frac{1}{z_j + \bar{z}_j} \right) + \frac{1}{2\pi i} \sum_{\substack{k=1 \\ k \neq j}}^N \Gamma_k \left( \frac{1}{z_j - z_k} - \frac{1}{z_j - \bar{z}_k} - \frac{1}{z_j + \bar{z}_k} + \frac{1}{z_j + z_k} \right)
\end{equation}


\medskip


\section{The $N$-vortex system in the upper right quadrant}
\label{sec:Nvor}

A vortex that remains uniformly separated from all other vortices
cannot approach either wall or the corner in finite time.

\subsection{Solitary Boundary Collisions and Energy Conservation}

For each vortex, define
\[
\rho_j=\frac{x_jy_j}{\sqrt{x_j^2+y_j^2}}.
\]
Since
\[
\frac{\min\{x_j,y_j\}}{\sqrt{2}}
\le \rho_j \le \min\{x_j,y_j\}
=\operatorname{dist}(z_j,\partial\Omega),
\]
approaching either wall or the corner is equivalent to
\(\rho_j\to0\).

The conserved Kirchhoff Hamiltonian \eqref{eq:1/4_Ham} is
\[
H=
\LaTeXunderbrace{\sum_{1\le j<k\le N}
\Gamma_j\Gamma_kG(z_j,z_k)}_{H_{\mathrm{int}}}
+
\LaTeXunderbrace{\frac{1}{8\pi}\sum_{j=1}^N
\Gamma_j^2\log(4\rho_j^2)}_{H_{\mathrm{self}}}.
\]
For \(\Gamma_j\ne0\), the individual self-energy
$
\frac{\Gamma_j^2}{8\pi}\log(4\rho_j^2)
$
tends to \(-\infty\) as vortex \(j\) approaches the boundary.
For a single vortex with nonzero circulation, conservation of \(H\)
therefore gives
$
\rho_1(t)=\rho_1(0)>0.
$
Consequently, a single vortex remains a positive distance from the
boundary, even as \(t\to\infty\). \cite{Yang}

For several vortices, an energy argument requires control of the
remaining terms. Suppose the configuration stays uniformly bounded,
\[
\sup_{0\le t<T}\max_{1\le k\le N}|z_k(t)|<\infty.
\]
Then all self-energy terms are bounded above. If a vortex with
nonzero circulation approaches the boundary along a sequence
\(t_n\uparrow T\), it follows that
\[
H_{\mathrm{self}}(t_n)\to-\infty
\implies
H_{\mathrm{int}}(t_n)\to+\infty.
\]
The quadrant Green's function satisfies
\[
G(z_a,z_b)=\frac{1}{4\pi}\log\!\left(
1+\frac{16x_ay_ax_by_b}
{|z_a-z_b|^2\,|z_a+z_b|^2}
\right).
\]
Thus \(G\ge0\), and \(G\) is uniformly bounded for pairs with a
positive lower bound on their separation in a bounded configuration.
After passing to a subsequence, there must therefore be a fixed
pair \(a\ne b\) such that
\[
\Gamma_a\Gamma_b>0,
\qquad
G(z_a(t_n),z_b(t_n))\to+\infty,
\qquad
|z_a(t_n)-z_b(t_n)|\to0.
\]


\subsection{Simultaneous Boundary and Inter-Vortex Collisions}

The energy argument leaves the close pair \((a,b)\) unspecified.
The following local estimate shows that the vortex approaching
the boundary must itself come arbitrarily close to another vortex.
When the trajectories have finite limits, this gives a simultaneous
collision at the same boundary point.


\begin{proposition}[(Finite-time approach to the boundary)]
Let \(z_1(t),\ldots,z_N(t)\) be a solution of
\eqref{eq:1/4_vortex_eq} with pairwise distinct positions in
\(\Omega\) for \(0\le t<T<\infty\). If
\[
\liminf_{t\uparrow T}
\operatorname{dist}(z_j(t),\partial\Omega)=0,
\]
then there exists \(k\ne j\) such that
\[
\liminf_{t\uparrow T}|z_j(t)-z_k(t)|=0.
\]
In particular, if all trajectories have finite limits at \(T\),
then vortex \(j\) and at least one other vortex converge to the
same boundary point.

\end{proposition}


\begin{proof}

Write \(r_j^2=x_j^2+y_j^2\). Let
\((U(x,y,t),V(x,y,t))\) denote the velocity field induced by all
vortices other than \(j\), together with their images, and set
\[
u_j=U(x_j,y_j,t),
\qquad
v_j=V(x_j,y_j,t).
\]
Differentiating \(\log\rho_j\) and substituting
\eqref{eq:1/4_vortex_eq} gives
\[
\begin{aligned}
\frac{d}{dt}\log\rho_j
&=\frac{y_j^2}{x_jr_j^2}\dot x_j
 +\frac{x_j^2}{y_jr_j^2}\dot y_j\\
&=\frac{\Gamma_j}{4\pi}
 \left(\frac{x_jy_j}{r_j^4}-\frac{x_jy_j}{r_j^4}\right)
 +\frac{y_j^2}{r_j^2}\frac{u_j}{x_j}
 +\frac{x_j^2}{r_j^2}\frac{v_j}{y_j}\\
&=\frac{y_j^2}{r_j^2}\frac{u_j}{x_j}
 +\frac{x_j^2}{r_j^2}\frac{v_j}{y_j}.
\end{aligned}
\]
Thus the self-image contributions cancel exactly.

Suppose there exist \(\delta>0\) and \(t_0<T\) such that
\[
|z_j(t)-z_k(t)|\ge\delta
\qquad
\text{for all }k\ne j,\quad t_0\le t<T.
\]
Every image of another vortex is also at least \(\delta\) from
\(z_j\), since reflection across either axis increases the squared
distance from a point in \(\Omega\). For example,
\[
|(x_j,y_j)-(-x_k,y_k)|^2
=|z_j-z_k|^2+4x_jx_k.
\]

If \(x_j\le\delta/2\), every point \(p=(s,y_j)\) with
\(0\le s\le x_j\) satisfies
\[
|p-a|
\ge |z_j-a|-|z_j-p|
\ge\delta-x_j
\ge\frac{\delta}{2}
\]
for every source \(a\) contributing to \((U,V)\).
The point-vortex kernel
\[
K(q)=\frac{(-q_2,q_1)}{2\pi|q|^2}
\]
satisfies
\[
|K(q)|=\frac{1}{2\pi|q|},
\qquad
\|DK(q)\|=\frac{1}{2\pi|q|^2},
\]
where the second norm is the Euclidean operator norm.
Hence \(\partial_xU\) is uniformly bounded on this segment.
The image construction gives \(U(0,y_j,t)=0\), so
\[
|u_j|
=\left|\int_0^{x_j}\partial_xU(s,y_j,t)\,ds\right|
\le Cx_j.
\]
If \(x_j>\delta/2\), separation from all sources gives a uniform
velocity bound \(|u_j|\le M\), and therefore
\[
|u_j|\le M\le\frac{2M}{\delta}x_j.
\]
The analogous argument along the vertical segment, using
\(V(x_j,0,t)=0\), gives \(|v_j|\le Cy_j\).
After increasing \(C\) if necessary, both estimates hold with a
constant depending only on \(\delta\) and the circulations.

Since \(x_j^2/r_j^2+y_j^2/r_j^2=1\), it follows that
\[
\left|\frac{d}{dt}\log\rho_j\right|\le C.
\]
Consequently,
\[
\rho_j(t)
\ge \rho_j(t_0)e^{-C(t-t_0)}
\ge \rho_j(t_0)e^{-C(T-t_0)}>0,
\qquad t_0\le t<T.
\]
The inequality
\(\operatorname{dist}(z_j,\partial\Omega)\ge\rho_j\)
then excludes a finite-time approach to the boundary.

Thus the assumed boundary approach precludes uniform separation
of \(z_j\) from all other vortices. Since there are finitely many
vortices, some fixed \(k\ne j\) satisfies
\[
\liminf_{t\uparrow T}|z_j(t)-z_k(t)|=0.
\]
If all trajectories have finite limits, this implies
\(z_j(T)=z_k(T)\in\partial\Omega\), proving the final assertion.

\end{proof}

These conclusions give necessary conditions for a boundary collision.
Establishing the existence of such a collision requires further
analysis.

\subsection{Internal Finite-Time Collisions ($N \ge 2$)}

\begin{theorem}[(No interior two-vortex collapse)]\label{thm:two_interior}
Let \(\Gamma_1,\Gamma_2\ne0\), and let \(z_1(t),z_2(t)\) solve
\eqref{eq:1/4_vortex_eq} with distinct positions in the quadrant
\(\Omega\) for \(0\le t<T<\infty\). Then, for every \(z^*\in\Omega\),
\[
\liminf_{t\uparrow T}
\left(|z_1(t)-z^*|+|z_2(t)-z^*|\right)>0.
\]
\end{theorem}

\begin{proof}
Suppose, for contradiction, that the limit inferior is zero for
some \(z^*\in\Omega\). There exists a sequence \(t_n\uparrow T\)
such that
\[
z_1(t_n)\to z^*,
\qquad
z_2(t_n)\to z^*.
\]
In particular,
\[
d_n:=|z_1(t_n)-z_2(t_n)|\to0.
\]

Write the Dirichlet Green's function as
\[
G(z,w)=-\frac{1}{2\pi}\log|z-w|+g(z,w),
\]
where, for the quadrant,
\[
g(z,w)=\frac{1}{2\pi}\log\!\left(
\frac{|z-\bar w|\,|z+\bar w|}{|z+w|}
\right).
\]
Since \(z^*\in\Omega\), the regular part \(g\) is smooth near
\((z^*,z^*)\), and the Robin function \(h(z)=g(z,z)\) is smooth
near \(z^*\).

With the convention in \eqref{eq:1/4_Ham}, conservation of the
Hamiltonian gives
\[
H(0)=H(t_n)
=-\frac{\Gamma_1\Gamma_2}{2\pi}\log d_n+R_n,
\]
where
\[
R_n=
\Gamma_1\Gamma_2g(z_1(t_n),z_2(t_n))
+\frac{\Gamma_1^2}{2}h(z_1(t_n))
+\frac{\Gamma_2^2}{2}h(z_2(t_n)).
\]
By continuity,
\[
R_n\longrightarrow
\frac{(\Gamma_1+\Gamma_2)^2}{2}h(z^*)\in\mathbb{R}.
\]
However, \(\log d_n\to-\infty\) and
\(\Gamma_1\Gamma_2\ne0\), so
\[
H(t_n)\longrightarrow
\begin{cases}
+\infty, & \Gamma_1\Gamma_2>0,\\
-\infty, & \Gamma_1\Gamma_2<0.
\end{cases}
\]
This contradicts \(H(t_n)=H(0)\in\mathbb{R}\), proving the claim.
\end{proof}

\paragraph{$N=3$ is the True Minimum.}
To achieve an internal collapse, the Hamiltonian requires a second logarithmic singularity to balance the first one, keeping the total energy finite.

If we add a third vortex, the interaction Hamiltonian contains:
$$H_{\text{int}} \propto \Gamma_1 \Gamma_2 \ln\vert{}z_1 - z_2\vert{} + \Gamma_2 \Gamma_3 \ln\vert{}z_2 - z_3\vert{} + \Gamma_3 \Gamma_1 \ln\vert{}z_3 - z_1\vert{}.$$
If all three vortices collapse toward $z^*$, one can select circulations of mixed signs (e.g., $\Gamma_1 \Gamma_2 > 0$ driving the energy to $-\infty$, while $\Gamma_2 \Gamma_3 < 0$ drives it to $+\infty$) and specific geometric configurations such that the infinities perfectly cancel, maintaining a constant $H$. For $N=2$, there is no third vortex available to provide this balancing term.


\medskip


\section{Finite-time Boundary Non-collision of Two vortices in the Quadrant}
\label{sec:BoundaryNoncollisionTwo}

\begin{theorem}[(No simultaneous boundary accumulation for two vortices)]
\label{thm:no-boundary-accumulation}
Let \(\Gamma_1,\Gamma_2\ne0\), and let \(z_1(t),z_2(t)\) solve
\eqref{eq:1/4_vortex_eq} with distinct positions in the quadrant
\(\Omega\) for \(0\le t<T<\infty\). Then, for every
\(p\in\partial\Omega\),
\[
\liminf_{t\uparrow T}
\left(|z_1(t)-p|+|z_2(t)-p|\right)>0.
\]
This includes points on either wall and the corner \(p=0\).
\end{theorem}

\begin{proof}
Suppose, for contradiction, that there exist \(p\in\partial\Omega\)
and a sequence \(t_n\uparrow T\) such that
\[
z_1(t_n)\to p,
\qquad
z_2(t_n)\to p.
\]
Set
\[
\rho_j=\frac{x_jy_j}{\sqrt{x_j^2+y_j^2}}.
\]
Since
\[
0<\rho_j\le\min\{x_j,y_j\}
=\operatorname{dist}(z_j,\partial\Omega),
\]
we have \(\rho_j(t_n)\to0\) for \(j=1,2\).

If \(\Gamma_1\Gamma_2<0\), positivity of the quadrant Green's
function and conservation of \eqref{eq:1/4_Ham} give
\[
H(0)=H(t_n)
\le \frac{1}{8\pi}\sum_{j=1}^{2}
\Gamma_j^2\log\bigl(4\rho_j(t_n)^2\bigr)
\longrightarrow-\infty.
\]
This contradicts \(H(0)\in\mathbb{R}\).

Now suppose \(\Gamma_1,\Gamma_2>0\), and define
\[
X(t)=\Gamma_1x_1(t)+\Gamma_2x_2(t),
\qquad
Y(t)=\Gamma_1y_1(t)+\Gamma_2y_2(t).
\]
Summing \eqref{eq:1/4_vortex_eq} and pairing the interaction terms
gives
\[
\begin{aligned}
\dot X
&=\frac{2}{\pi}\sum_{j,k=1}^{2}
\frac{\Gamma_j\Gamma_kx_jx_k(y_j+y_k)}
{d_{2,jk}^{2}d_{4,jk}^{2}}>0,\\
\dot Y
&=-\frac{2}{\pi}\sum_{j,k=1}^{2}
\frac{\Gamma_j\Gamma_ky_jy_k(x_j+x_k)}
{d_{3,jk}^{2}d_{4,jk}^{2}}<0.
\end{aligned}
\]
The terms with \(j=k\) account for the self-image contributions.

If \(p=(0,b)\) with \(b\ge0\), including the corner, then
\(X(t_n)\to0\), contradicting \(X(t)\ge X(0)>0\).
Thus the only remaining possibility is \(p=(a,0)\) with \(a>0\).

Along the selected sequence,
\[
X(t_n)\to(\Gamma_1+\Gamma_2)a,
\qquad
Y(t_n)\to0.
\]
By monotonicity, these are also the full limits,
\[
X(t)\longrightarrow(\Gamma_1+\Gamma_2)a,
\qquad
Y(t)\longrightarrow0
\qquad\text{as }t\uparrow T.
\]
Because the circulations are positive, the positions remain
bounded and \(y_1(t),y_2(t)\to0\). Consequently,
\[
H_{\mathrm{self}}(t)\to-\infty.
\]
Conservation of \(H\) and \(\Gamma_1\Gamma_2>0\) therefore imply
\[
G(z_1(t),z_2(t))\to+\infty.
\]

The identity
\[
G(z_1,z_2)=\frac{1}{4\pi}\log\!\left(
1+\frac{16x_1y_1x_2y_2}
{|z_1-z_2|^2\,|z_1+z_2|^2}
\right)
\]
shows that, for bounded positions, \(G\) is uniformly bounded
whenever \(|z_1-z_2|\) is bounded away from zero.
Here we use \(|z_1+z_2|\ge|z_1-z_2|\) in the quadrant.
It follows that
\[
|z_1(t)-z_2(t)|\to0.
\]
Combining this with the limit of \(X(t)\) yields
\[
x_1(t)\to a,
\qquad
x_2(t)\to a.
\]

Choose \(t_0<T\) and \(c>0\) such that \(x_j(t)\ge c\) for
\(j=1,2\) and \(t_0\le t<T\). Since
\[
d_{3,jk}^{2},d_{4,jk}^{2}\ge(x_j+x_k)^2,
\]
the identity for \(\dot Y\) gives
\[
0\le-\dot Y
\le\frac{2}{\pi(2c)^3}
\sum_{j,k=1}^{2}\Gamma_j\Gamma_ky_jy_k
=\frac{Y^2}{4\pi c^3}.
\]
Writing \(C=(4\pi c^3)^{-1}\) and integrating, we obtain
\[
Y(t)\ge
\frac{Y(t_0)}{1+C\,Y(t_0)(t-t_0)}
\ge
\frac{Y(t_0)}{1+C\,Y(t_0)(T-t_0)}>0.
\]
This contradicts \(Y(t)\to0\).

Finally, if \(\Gamma_1,\Gamma_2<0\), the transformation
\[
\widetilde z_j=i\overline{z_j}=y_j+ix_j,
\qquad
\widetilde\Gamma_j=-\Gamma_j>0
\]
maps the system to the same quadrant equations with positive
circulations. It also preserves distances and maps
\(\partial\Omega\) to itself, reducing this case to the preceding
argument. The contradiction proves the claim.
\end{proof}

\begin{remark}
We have proved the following monotonicity property above: 
for positive vortices,
$$X(t)\ \text{is strictly increasing},\qquad Y(t)\ \text{is strictly decreasing}.$$
This is precisely the point-vortex counterpart of the first-moment monotonicity \(X_1'(t)\ge0\), \(X_2'(t)\le0\) for nonnegative odd-odd Euler vorticity in the quadrant.
\cite{choi2024stability}
\end{remark}

By Theorem \ref{thm:no-boundary-accumulation}, and ruling out escape to infinity in finite time, we can prove the global existence holds for two vortices with arbitrary nonzero circulations and distinct initial positions in the quadrant.

\begin{theorem}[(Global existence for two vortices in the quadrant)]
\label{thm:two_global}
Let \(\Gamma_1,\Gamma_2\in\mathbb{R}\setminus\{0\}\), and let
\(z_1^0,z_2^0\in\Omega\) be distinct. The solution of
\eqref{eq:1/4_vortex_eq} with \(z_j(0)=z_j^0\) exists uniquely
for all \(t\in\mathbb{R}\).

Moreover, for every \(S>0\),
\[
\inf_{|t|\le S}
\min\left\{
x_1(t),y_1(t),x_2(t),y_2(t),|z_1(t)-z_2(t)|
\right\}>0,
\]
and
\[
\sup_{|t|\le S}\max_{j=1,2}|z_j(t)|<\infty.
\]
\end{theorem}

\begin{proof}
Local existence and uniqueness hold on
$
\mathcal C=\{(z_1,z_2)\in\Omega^2:z_1\ne z_2\}.
$
Suppose the maximal forward existence time is \(T<\infty\).
Write
\[
d(t)=|z_1(t)-z_2(t)|,
\qquad
\rho_j(t)=\frac{x_j(t)y_j(t)}
{\sqrt{x_j(t)^2+y_j(t)^2}}.
\]
We first establish a uniform positive lower bound on all four
coordinates for \(0\le t<T\).

For like-sign circulations, the reflection symmetry used earlier
allows us to assume \(\Gamma_1,\Gamma_2>0\). The weighted-coordinate
identity for
$
Y=\Gamma_1y_1+\Gamma_2y_2
$
gives \(\dot Y<0\), so \(y_j\le M\) for some constant \(M\).
Consequently, the self-energy terms are bounded above and
\[
G(z_1,z_2)\le
\frac{1}{4\pi}\log\!\left(1+\frac{4M^2}{d^2}\right).
\]
Thus, along any sequence on which some \(\rho_j\to0\),
energy conservation forces \(G\to+\infty\) and hence \(d\to0\).

If some \(x_j(t_n)\to0\), both horizontal coordinates tend to zero.
Since the vertical coordinates are bounded, a subsequence gives
a common boundary accumulation point, contradicting
Theorem~\ref{thm:no-boundary-accumulation}. Therefore \(x_j\ge c>0\).
The weighted-coordinate estimate from the earlier proof now gives
\[
\dot Y\ge-\frac{Y^2}{4\pi c^3},
\qquad
Y(t)\ge
\frac{Y(0)}{1+Y(0)T/(4\pi c^3)}>0.
\]
If some \(y_j(t_n)\to0\), the same energy argument forces
\(d(t_n)\to0\), and hence \(Y(t_n)\to0\), a contradiction.
All coordinates therefore have positive lower bounds.

For opposite-sign circulations, \(G\ge0\) and conservation of \(H\)
give
\[
\max\{\rho_1,\rho_2\}\ge
c_0:=
\frac12\exp\!\left(
\frac{4\pi H}{\Gamma_1^2+\Gamma_2^2}
\right)>0.
\]
Set \(b=c_0/(2\sqrt2)\). Whenever \(\rho_j\le b\), the inequalities
\[
\rho_k\le\min\{x_k,y_k\},
\qquad
\min\{x_j,y_j\}\le\sqrt2\,\rho_j,
\qquad k\ne j,
\]
imply \(d\ge c_0/2\). The earlier local estimate
\[
\left|\frac{d}{dt}\log\rho_j\right|
\le\frac{C_\Gamma}{d^2}
\]
therefore bounds this derivative whenever \(\rho_j\le b\).
Integrating on each interval below that threshold yields
\[
\rho_j(t)\ge
\min\{\rho_j(0),b\}e^{-Ct}
\ge\min\{\rho_j(0),b\}e^{-CT}>0.
\]
Since \(x_j,y_j\ge\rho_j\), we have established, in every case,
\[
x_j(t),y_j(t)\ge\eta>0
\qquad(0\le t<T,\ j=1,2).
\]

Let \(F(q,t)\) denote the velocity field generated by all six image
vortices. Every point on the segment joining \(z_1(t)\) and
\(z_2(t)\) is at least \(2\eta\) from every image. The Biot--Savart
kernel bounds therefore give
\[
|F(z_j,t)|\le C,
\qquad
|F(z_1,t)-F(z_2,t)|\le Ld,
\]
uniformly for \(t<T\), independently of the sizes of the positions.

Using real-vector notation, the direct mutual interaction is
perpendicular to \(z_1-z_2\). Hence
\[
\frac{d}{dt}d^2
=2(z_1-z_2)\cdot
\bigl(F(z_1,t)-F(z_2,t)\bigr),
\]
so \(|\dot d|\le Ld\) and
\[
d(t)\ge d(0)e^{-LT}=:\delta>0.
\]
The full vortex velocities are consequently bounded,
\[
|\dot z_j(t)|
\le\frac{|\Gamma_k|}{2\pi\delta}+C
\le V,
\qquad k\ne j.
\]
In particular,
\[
|z_j(t)|\le|z_j(0)|+VT,
\qquad 0\le t<T,
\]
which excludes finite-time escape to infinity.

These bounds confine the trajectory to a compact subset of
\(\mathcal C\). The ODE continuation theorem extends the solution
beyond \(T\), contradicting maximality. Reversing time changes both
circulation signs, so the same argument gives backward global
existence.
\end{proof}

\medskip

\section{Finite-time boundary exclusion for three vortices}
\label{sec:three-boundary-exclusion}

We consider classical solutions with nonzero circulations and
pairwise distinct positions in
\[
\Omega=\{z=x+iy:x>0,\ y>0\},
\qquad 0\le t<T<\infty.
\]
No boundedness or trajectory-convergence assumption is imposed.
Constants denoted by \(C\) may change from line to line.

\subsection{Basic identities}

Identifying complex numbers with vectors in \(\mathbb R^2\), set
\[
K(z)=\frac{iz}{2\pi|z|^2},
\qquad
B(z,w)=-K(z-\bar w)-K(z+\bar w)+K(z+w),
\qquad
\mathcal K(z,w)=K(z-w)+B(z,w).
\]
The quadrant equations are
\begin{equation}
\dot z_j
=
\Gamma_jB(z_j,z_j)
+\sum_{k\ne j}\Gamma_k\mathcal K(z_j,z_k).
\label{eq:selfcontained-quadrant}
\end{equation}
They conserve
\begin{equation}
H=\sum_{j<k}\Gamma_j\Gamma_kG(z_j,z_k)
+\frac1{4\pi}\sum_j\Gamma_j^2\log(2\rho_j),
\qquad
\rho_j=\rho(z_j),\quad
\rho(z)=\frac{xy}{|z|}.
\label{eq:selfcontained-H}
\end{equation}
For \(z=x+iy\), \(w=u+iv\),
\[
G(z,w)=\frac1{4\pi}
\log\left(1+\frac{16xyuv}{|z-w|^2|z+w|^2}\right).
\]
In particular,
\begin{equation}
0\le G(z,w)\le
\frac1{4\pi}
\log\left(
1+\frac{\min\{4xu,4yv,8\rho(z)\rho(w)\}}{|z-w|^2}
\right).
\label{eq:selfcontained-G-bound}
\end{equation}
Also,
\begin{equation}
\rho\le\min\{x,y\}\le\sqrt2\,\rho,\qquad
2\rho\le|z|,\qquad
\nabla\rho=\frac{(y^3,x^3)}{|z|^3},\qquad
|\nabla\rho|\le1,\qquad
\|D^2\rho\|=\frac{3\rho}{|z|^2}.
\label{eq:selfcontained-rho}
\end{equation}
Matrix norms above denote the operator norm induced by the Euclidean norm

\begin{lemma}[(Normal-field estimates and scalar barriers)]
\label{lem:summary-local-estimates}
Suppose that a vortex is at least a distance \(\delta>0\) from
each vortex in a specified collection. If \((u,v)\) is the
velocity induced by that collection and its images, then
\begin{equation}
|u|\le C_\delta x,
\qquad
|v|\le C_\delta y.
\label{eq:selfcontained-normal}
\end{equation}
If \((u_j,v_j)\) denotes the velocity induced by all vortices
other than \(j\) and their images, then
\begin{equation}
\frac{d}{dt}\log\rho(z_j)
=
\frac{y_j^2}{x_j^2+y_j^2}\frac{u_j}{x_j}
+
\frac{x_j^2}{x_j^2+y_j^2}\frac{v_j}{y_j}.
\label{eq:selfcontained-rho-derivative}
\end{equation}

We shall also use the following scalar facts.
If \(f>0\) is locally absolutely continuous and
\[
f'\ge-Cf
\quad\text{almost everywhere whenever }f\le\varepsilon,
\]
then
\begin{equation}
f(t)\ge \min\{f(t_0),\varepsilon\}e^{-C(t-t_0)}
\qquad(t_0\le t<T).
\label{eq:summary-scalar-barrier}
\end{equation}
Furthermore, suppose that \(\liminf_{t\uparrow T}f(t)=0\).
If, on every sufficiently small fixed strip
\(\varepsilon\le f\le2\varepsilon\), either
\(f'\le C_\varepsilon\) or \(f'\ge-C_\varepsilon\) holds
sufficiently late in time, then \(f(t)\to0\).
\end{lemma}

\begin{proof}
Every image of a source is at least as far from a point in
\(\Omega\) as the source itself. If \(x\le\delta/2\), the
horizontal segment to the vertical wall stays at least
\(\delta/2\) from all sources and their images. The induced
field has bounded spatial derivatives there, and its horizontal
component vanishes on the wall. The mean value theorem gives
\(|u|\le C_\delta x\). If \(x>\delta/2\), the same conclusion
follows from the uniform velocity bound at the target point.
The vertical estimate is identical.

Differentiating \(\log\rho\) gives
\[
\frac{d}{dt}\log\rho
=
\frac{y^2}{x^2+y^2}\frac{\dot x}{x}
+
\frac{x^2}{x^2+y^2}\frac{\dot y}{y}.
\]
The self-image velocity is
\[
\frac{\Gamma}{4\pi}
\left(
\frac{x^2}{y(x^2+y^2)},
-\frac{y^2}{x(x^2+y^2)}
\right).
\]
Its contribution to this derivative is
\[
\frac{\Gamma xy}{4\pi(x^2+y^2)^2}
-
\frac{\Gamma xy}{4\pi(x^2+y^2)^2}
=0,
\]
proving \eqref{eq:selfcontained-rho-derivative}.

The scalar barrier follows by integration on each connected
component of \(\{f<\varepsilon\}\). For the last assertion, an
upward crossing of \([\varepsilon,2\varepsilon]\) under the
upper derivative bound, or a downward crossing under the lower
derivative bound, takes at least
\(\varepsilon/C_\varepsilon\) units of time. Thus such crossings
cannot occur arbitrarily close to the finite time \(T\).
Together with \(\liminf f=0\), this proves \(f\to0\).
\end{proof}

For a group \(I\), define
\(X_I=\sum_{j\in I}\Gamma_jx_j\) and
\(Y_I=\sum_{j\in I}\Gamma_jy_j\).
Let \((U_j,V_j)\) denote the velocity induced at \(z_j\) by the vortices outside \(I\) and their images (additional velocity), pairing
the equations yields
\begin{align}
X_I'
&=
\frac2\pi\sum_{j,k\in I}
\frac{\Gamma_j\Gamma_kx_jx_k(y_j+y_k)}
 {|z_j-\bar z_k|^2|z_j+z_k|^2}
+\sum_{j\in I}\Gamma_jU_j,
\label{eq:selfcontained-X}\\
Y_I'
&=
-\frac2\pi\sum_{j,k\in I}
\frac{\Gamma_j\Gamma_ky_jy_k(x_j+x_k)}
 {|z_j+\bar z_k|^2|z_j+z_k|^2}
+\sum_{j\in I}\Gamma_jV_j.
\label{eq:selfcontained-Y}
\end{align}
The sums include self-images through \(j=k\).
For the full system, the elementary Laplace transform of a
product of sines gives
\begin{align}
\left(\sum_j\Gamma_jx_j\right)'
&=
\frac1\pi\int_0^\infty
\left(\sum_j\Gamma_j e^{-s y_j}\sin(sx_j)\right)^2ds
\ge0,
\label{eq:selfcontained-signed-X}\\
\left(\sum_j\Gamma_jy_j\right)'
&=
-\frac1\pi\int_0^\infty
\left(\sum_j\Gamma_j e^{-s x_j}\sin(sy_j)\right)^2ds
\le0.
\label{eq:selfcontained-signed-Y}
\end{align}

\subsection{Auxiliary lemmas}

\begin{lemma}[(A forced same-sign group)]
\label{lem:three-forced-group}
Let \(I\) be a finite nonempty group of vortices with circulations
of one sign. Suppose that their motion consists of their mutual
quadrant vortex dynamics and additional velocities \((U_j,V_j)\).
Assume that, for some \(b,L>0\),
\[
|U_j|\le Lx_j\quad\text{when }x_j\le b,
\qquad
|V_j|\le Ly_j\quad\text{when }y_j\le b.
\]
Then all its coordinates have positive
lower bounds on \([0,T)\), i.e.,
\[
\inf_{0\le t<T}\min_{j\in I}
\operatorname{dist}(z_j(t),\partial\Omega)>0.
\]
\end{lemma}

\begin{proof}
The symmetry
\((x,y,\Gamma)\mapsto(y,x,-\Gamma)\) reduces the proof to
positive circulations \(a_j\).

For every nonempty subset \(Q\) of the group, put
\(X_Q=\sum_{j\in Q}a_jx_j\).
We first establish positive lower bounds for all \(X_Q\), by
descending induction on \(|Q|\). Put
\[
m=\min_{j\in I}a_j,\qquad M=\max_{j\in I}a_j.
\]
When \(X_I\le mb\), every \(x_j\le b\), and
\eqref{eq:selfcontained-X} gives
\(\dot X_I\ge-LX_I\). The scalar barrier therefore bounds
\(X_I\) away from zero.

Suppose that every subset of size \(q+1\) has its horizontal
weighted sum bounded below by \(A_0>0\). For \(|Q|=q\), set
\[
\delta=\min\left\{b,\frac{A_0}{4M}\right\},
\qquad \varepsilon=m\delta.
\]
Whenever \(X_Q\le\varepsilon\),
\[
x_j\le\delta\quad(j\in Q),
\qquad
x_k\ge3\delta\quad(k\in I\setminus Q).
\]
Indeed,
\(a_kx_k=X_{Q\cup\{k\}}-X_Q\ge3A_0/4\).
The omitted vortices therefore induce
\(|u_j^Q|\le C_\delta x_j\), by the wall-segment argument in
Lemma~\ref{lem:summary-local-estimates}. Consequently,
\[
\dot X_Q\ge-CX_Q\qquad(X_Q\le\varepsilon).
\]
The scalar barrier completes the induction. In particular,
\(x_j\ge c_0>0\) for all \(j\in I\).

With this lower bound,
\[
0\le
\frac2\pi\sum_{j,k\in Q}
\frac{a_ja_ky_jy_k(x_j+x_k)}
     {d_{3,jk}^2d_{4,jk}^2}
\le
\frac{Y_Q^2}{4\pi c_0^3}.
\]
We now repeat the descending induction for \(Y_Q\).
For the full group, sufficiently small \(Y_I\) gives
\[
\dot Y_I\ge-LY_I-\frac{Y_I^2}{4\pi c_0^3}\ge-CY_I.
\]
At the induction step, small \(Y_Q\) separates the members of
\(Q\) vertically from all omitted vortices. Their normal field
satisfies \(|v_j^Q|\le Cy_j\), and hence
\[
\dot Y_Q\ge-CY_Q-\frac{Y_Q^2}{4\pi c_0^3}
\ge-C'Y_Q
\]
below a fixed threshold. The scalar barrier gives positive lower
bounds for every \(Y_Q\), and therefore for every \(y_j\).
\end{proof}

\begin{lemma}[(A cutoff at infinity)]
\label{lem:three-cutoff}
There is no sequence \(t_n\uparrow T\) such that
\[
\rho_1(t_n),\rho_2(t_n)\to0,
\qquad
\rho_3(t_n)\to\infty.
\]
The same assertion holds after any relabeling.
\end{lemma}

\begin{proof}
Choose a smooth \(\chi:\mathbb R\to[0,1]\), equal to zero
for \(s\le1\) and one for \(s\ge2\), and set
\[
\phi_R(x,y)=\chi(x/R)\chi(y/R),
\qquad
M_R=\sum_j\Gamma_j\phi_R(z_j).
\]
Then \(\|\nabla\phi_R\|_\infty\le C/R\) and
\(\|D^2\phi_R\|_\infty\le C/R^2\).
Pairing the direct interactions gives
\[
\begin{aligned}
M_R'
={}&
\sum_{j<k}\Gamma_j\Gamma_k
\bigl(\nabla\phi_R(z_j)-\nabla\phi_R(z_k)\bigr)
\cdot K(z_j-z_k)\\
&+
\sum_{j,k}\Gamma_j\Gamma_k
\nabla\phi_R(z_j)\cdot B(z_j,z_k).
\end{aligned}
\]
By the mean value theorem and the Hessian bound,
\[
|\nabla\phi_R(z_j)-\nabla\phi_R(z_k)|
\le \frac{C}{R^2}|z_j-z_k|.
\]
Since \(|K(z_j-z_k)|=(2\pi|z_j-z_k|)^{-1}\), it follows that
\[
\left|
\bigl(\nabla\phi_R(z_j)-\nabla\phi_R(z_k)\bigr)
\cdot K(z_j-z_k)
\right|
\le \frac{C}{R^2}.
\]
Thus the first sum is \(O(R^{-2})\), independently of how
small the inter-vortex distances become.

For the second sum, \(\chi=\chi'=0\) on \((-\infty,1]\)
implies that \(\nabla\phi_R(z_j)\ne0\) only when
\(x_j,y_j>R\). Every image source is then at least distance
\(R\) from \(z_j\), including self-images, because
\[
|z_j-\bar z_k|\ge y_j+y_k>R,\qquad
|z_j+\bar z_k|,\ |z_j+z_k|\ge x_j+x_k>R.
\]
Hence \(|B(z_j,z_k)|\le3/(2\pi R)\). Together with
\(|\nabla\phi_R|\le C/R\), this bounds each nonzero image
summand by \(C|\Gamma_j\Gamma_k|/R^2\).
Summing both contributions gives
\[
|M_R'(t)|\le \frac{C_\Gamma}{R^2},
\]
where \(C_\Gamma\) depends only on the cutoff and the
circulations, not on the vortex positions or separations.

Choose \(R\) large enough that \(M_R(0)=0\) and
\(C_\Gamma T/R^2<|\Gamma_3|\). Then
\[
|M_R(t)|<|\Gamma_3|.
\]
Along the proposed sequence, however,
\(\phi_R(z_1)=\phi_R(z_2)=0\) and \(\phi_R(z_3)=1\)
eventually, by \eqref{eq:selfcontained-rho}.
Thus \(M_R(t_n)=\Gamma_3\), a contradiction.
\end{proof}

\begin{remark}
Think of \(\phi_R\) as a smooth indicator of the region far from both walls:
\[
\phi_R(x,y)=
\begin{cases}
0,&\min\{x,y\}\le R,\\
1,&\min\{x,y\}\ge2R,
\end{cases}
\qquad 0\le\phi_R\le1.
\]
Thus,
$
M_R(t)=\sum_j\Gamma_j\phi_R(z_j(t))
$
measures the signed circulation in that region, with fractional weights in the transition zone. A vortex near either wall contributes zero; a vortex sufficiently far from both walls contributes its full circulation. Because the circulations can have different signs, \(M_R\) need not be positive.

The important estimate is
\[
|M_R'(t)|\le \frac{C_\Gamma}{R^2}.
\]
It says this weighted circulation changes slowly when \(R\) is large. Pairwise cancellation makes the estimate valid even during close vortex encounters.
\end{remark}

For mixed signs, reflection and relabeling allow us to write
\[
(\Gamma_1,\Gamma_2,\Gamma_3)=(a,b,-c),
\qquad a,b,c>0.
\]
For the rest of the section, use
\[
A=a+b,\qquad Z=X+iY=\frac{az_1+bz_2}{A},
\qquad r=|z_1-z_2|,\qquad \ell=|Z-z_3|,
\]
\[
h=\rho(Z),\qquad s=\rho_3,\qquad
M=\max\{h,s\},\qquad d=\min\{h,\ell\}.
\]

\begin{lemma}[(Uniform concentration and effective-pair estimates)]
\label{lem:three-effective}
For sufficiently small \(M\),
\begin{equation}
r/h\le CM^\vartheta,\qquad
r\le Ch^\mu s^\nu,
\label{eq:three-concentration}
\end{equation}
where
\[
\vartheta=\frac{a^2+b^2+c^2}{2ab}>0,\qquad
\mu=\frac{A^2}{2ab}\ge2,\qquad
\nu=\frac{c^2}{2ab}>0,\qquad
\kappa=\mu+\nu>2.
\]
Moreover, \(\max_j\rho_j(t_n)\to0\) implies \(M(t_n)\to0\).

If \(r\le d/4\), replacing the positive pair by a vortex
of circulation \(A\) at \(Z\) gives
\[
\dot Z=V_Z+\mathcal R_Z,\qquad
\dot z_3=V_3+\mathcal R_3,
\]
where \(V_Z,V_3\) are the quadrant velocities of vortices
\(A,-c\) at \(Z,z_3\), and
\begin{equation}
|\mathcal R_Z|\le Cr^2/d^3,\qquad
|\mathcal R_3|\le Cr^2/\ell^3,\qquad
\frac{|\mathcal R_{3,x}|}{x_3}
+\frac{|\mathcal R_{3,y}|}{y_3}
\le Cr^2/\ell^4.
\label{eq:three-remainders}
\end{equation}
Consequently,
\begin{equation}
\mathcal E=
-AcG(Z,z_3)
+\frac{A^2}{4\pi}\log(2h)
+\frac{c^2}{4\pi}\log(2s)
\quad\text{satisfies}\quad
|\mathcal E'|\le Cr^2/d^4.
\label{eq:three-effective-energy}
\end{equation}
If \(c=A\), sufficiently small \(M\) also implies
\begin{equation}
d>0,\qquad r\le Cd^2,\qquad r\le d/4.
\label{eq:three-neutral}
\end{equation}
All constants are independent of spatial bounds.
\end{lemma}

\begin{proof}
Since \(\rho\) is coordinatewise increasing and homogeneous,
positivity of the centroid weights gives
\(\rho_1,\rho_2\le Ch\).
Discarding the nonpositive interactions from \(H=E\), we get
\[
4\pi E\le
ab\log(1+Ch^2/r^2)
+(a^2+b^2+c^2)\log(CM).
\]
Exponentiation proves the first estimate in
\eqref{eq:three-concentration}.

For \(r\le h/4\), use
\[
G(z,w)=-\frac1{2\pi}\log|z-w|+g(z,w),
\qquad
g(z,w)=\frac1{2\pi}
\log\frac{|z-\bar w|\,|z+\bar w|}{|z+w|}.
\]
Expansion of the regular terms about \(Z\) has uniformly
bounded error, since the image distances are bounded below
by a fixed multiple of \(h\). Therefore
\[
4\pi E\le
-2ab\log r+A^2\log(2h)+c^2\log(2s)+C,
\]
which proves the second estimate.
If \(m=\max_j\rho_j\to0\), the same argument using
\(\rho_1\rho_2\le m^2\) gives
\(r\le Cm^{1+\vartheta}\to0\).
Since \(\rho\) is \(1\)-Lipschitz,
\(h\le\rho_1+r\to0\), proving the additional assertion.

The direct pair velocities cancel in the centroid, and
\[
a(z_1-Z)+b(z_2-Z)=0
\]
cancels first-order Taylor terms.
The bounds \(|D^mK(w)|\le C_m|w|^{-m-1}\) give the ordinary
remainders in \eqref{eq:three-remainders}.
If \(x_3\le\ell/4\), the segment to the vertical wall stays
separated from the interpolated sources. The normal remainder
vanishes at the wall; one target derivative gives
\[
|\mathcal R_{3,x}|\le Cx_3r^2/\ell^4.
\]
If \(x_3>\ell/4\), use the ordinary remainder.
The vertical argument is identical.

The effective velocities conserve \(\mathcal E\).
Its gradient at \(Z\) is \(O(d^{-1})\).
At \(z_3\), the interaction gradient is \(O(\ell^{-1})\),
and the self-energy derivatives are bounded by
\(C/x_3,C/y_3\). This proves
\eqref{eq:three-effective-energy}.

Finally, let \(c=A\).
Equation \eqref{eq:three-concentration} gives \(r\le Ch^2\).
If \(\ell\ge h/4\), this proves \eqref{eq:three-neutral}.
If \(\ell<h/4\), all three positions lie within a fixed
fraction of \(h\) of \(Z\).
The common regular-energy coefficient is
\(\frac12(a+b-A)^2=0\), so
\[
2\pi E\le-ab\log r+A^2\log(\ell+r)+C.
\]
Hence \(r\le C(\ell+r)^{A^2/(ab)}\).
Since \(A^2/(ab)\ge4\), smallness forces \(\ell>r\), then
\(r\le C\ell^{A^2/(ab)}\).
Together with \(r\le Ch^2\), this proves the claim.
\end{proof}

\begin{lemma}[(No common corner accumulation)]
\label{lem:three-corner}
For a three-vortex solution,
\[
\liminf_{t\uparrow T}\max_j|z_j(t)|>0.
\]
\end{lemma}

\begin{proof}
The same-sign case follows from
Lemma~\ref{lem:three-forced-group}.
Use the mixed-sign notation above and suppose that
\(R(t_n):=\max_j|z_j(t_n)|\to0\).

By \eqref{eq:selfcontained-signed-X}--%
\eqref{eq:selfcontained-signed-Y},
\[
P=AX-cx_3\uparrow0,\qquad Q=AY-cy_3\downarrow0.
\]
Put \(\alpha=A/c\) and \(f=X+y_3\).
The relations
\[
x_3=\alpha X-P/c,\qquad Y=y_3/\alpha+Q/A
\]
give
\[
f\le2R,\qquad
R\le C(f+|P|+|Q|),\qquad f\le C|Z|.
\]
We show the derivation of the last two inequalities briefly. Positivity of the centroid weights gives
\(|z_1|\le(A/a)|Z|\) and \(|z_2|\le(A/b)|Z|\).
Thus, with \(C_0=\max\{A/a,A/b\}\),
\[
R\le\max\{C_0(X+Y),\,x_3+y_3\}.
\]
Since
\[
X+Y\le\max\{1,\alpha^{-1}\}f+\frac{|Q|}{A},
\qquad
x_3+y_3\le\max\{1,\alpha\}f+\frac{|P|}{c},
\]
we obtain \(R\le C(f+|P|+|Q|)\).
Moreover, \(Q\ge0\) implies \(y_3\le\alpha Y\), so
\[
f=X+y_3\le X+\alpha Y
\le\sqrt{1+\alpha^2}\,|Z|.
\]
Here \(C\) depends only on the circulations.

We show that \(f'\ge-C_\varepsilon\) on every sufficiently
small fixed strip \(\varepsilon\le f\le2\varepsilon\),
sufficiently late in time. There \(R\le C\varepsilon\).

If \(A\ne c\), then
\[
Z-z_3=(1-\alpha)Z+(P+iQ)/c
\]
implies \(\ell\ge c_\varepsilon>0\), with
\(c_\varepsilon\) a fixed multiple of \(\varepsilon\).
Concentration gives \(r\le CR^\kappa\ll\varepsilon\).
Thus both opposite-sign fields are bounded on the strip.
The positive pair's internal contribution to \(X'\), and
the negative vortex's self-image contribution to \(y_3'\),
are nonnegative. Hence \(f'\ge-C_\varepsilon\).

If \(A=c\), then
\[
u=x_3-X\ge0,\qquad v=Y-y_3\ge0.
\]
The neutral estimate gives \(r\le Cd^2\), so the effective
velocity errors are \(O(d)=O(R)\).
When \(X\ge\varepsilon/2\), the extra vertical-wall images
are bounded. For the upper-half-plane field, with
\(D_\pm=u^2+(Y\pm y_3)^2\),
\[
V^{\mathbb H}_{3,y}
=\frac{Au}{2\pi}(D_-^{-1}-D_+^{-1})\ge0,
\]
\[
V^{\mathbb H}_{Z,x}
=\frac{A}{2\pi}
\left(
\frac1{2Y}+\frac{Y-y_3}{D_-}
-\frac{Y+y_3}{D_+}
\right)\ge0.
\]
The last inequality follows by clearing denominators and
using \(Y\ge y_3\).
When \(y_3\ge\varepsilon/2\), the corresponding
right-half-plane formulas are obtained by replacing
\((Y,y_3,u)\) with \((x_3,X,v)\) and interchanging the
components. They give the same lower bound for \(f'\).

The crossing observation now gives \(f\to0\), and hence
\(R\to0\).
Indeed, since \(f(t_n)\le2R(t_n)\to0\), we have
\(\liminf_{t\uparrow T}f(t)=0\). If \(f\not\to0\), choose
a sufficiently small \(\varepsilon>0\) such that
\(\limsup_{t\uparrow T}f(t)>2\varepsilon\).
By continuity, there exist intervals \([a_n,b_n]\), with
\(a_n\uparrow T\), on which
\[
f(a_n)=2\varepsilon,\qquad f(b_n)=\varepsilon,\qquad
\varepsilon\le f(t)\le2\varepsilon.
\]
For sufficiently large \(n\), integrating
\(f'\ge-C_\varepsilon\) over these intervals gives
\[
\varepsilon
=f(a_n)-f(b_n)
\le C_\varepsilon(b_n-a_n)
\le C_\varepsilon(T-a_n)\longrightarrow0,
\]
a contradiction. Therefore \(f(t)\to0\) as \(t\uparrow T\).

If \(A\ne c\), the inequalities
\(x_3\ge\alpha X\), \(y_3\le\alpha Y\) imply
\[
\ell\ge
\begin{cases}
(\alpha-1)X,&\alpha>1,\\
(1-\alpha)Y,&\alpha<1.
\end{cases}
\]
Thus \(\ell\ge c_0h\), and concentration gives
\(r\le Cd^2\), \(r\le d/4\) near \(T\).
For \(A=c\), these are precisely
\eqref{eq:three-neutral}.
Therefore \(|\mathcal E'|\le C\), whereas
\[
\mathcal E\le
\frac{A^2}{4\pi}\log(2h)
+\frac{c^2}{4\pi}\log(2s)\longrightarrow-\infty.
\]
This is impossible in finite time.
\end{proof}

\begin{lemma}[(A modified weighted distance)]
\label{lem:three-modified-moment}
For the mixed-sign system, there is a continuously
differentiable function \(V(t)\) such that, whenever
\(M\) is sufficiently small,
\begin{equation}
|V'|\le CM,\qquad V=Ah-cs+o(M),
\label{eq:three-modified-moment}
\end{equation}
with a uniform \(o(M)\).
\end{lemma}

\begin{proof}
Write
\[
L=Ah-cs,\qquad W=a\rho_1+b\rho_2-cs.
\]
Lemma~\ref{lem:three-corner} supplies
\(R\ge r_0>0\) sufficiently late.
Since
\[
|z_1|,|z_2|\le C_*|Z|,
\qquad C_*=\frac{A}{\min(a,b)},
\]
we may choose \(\delta>0\) sufficiently small that
\[
|Z|\le2\delta\quad\Longrightarrow\quad\ell\ge r_0/2.
\]
Let \(\theta(Z)\) be smooth, equal to \(1\) on
\(|Z|\le\delta\), and zero on \(|Z|\ge2\delta\). Set
\[
V=\theta L+(1-\theta)W.
\]
Weighted Taylor expansion and concentration give
\[
|L-W|\le Cr^2/h
\le CM^{2\kappa-1}=o(M).
\]

Where \(|Z|\le2\delta\), the third vortex is separated.
The normal-field and effective-pair estimates give
\[
|h'|\le Ch+Cr^2/h^3\le Ch,\qquad |s'|\le Cs,
\qquad |\dot Z|\le C/h.
\]
Here the estimate for \(s'\) uses the normal remainder in
\eqref{eq:three-remainders}. Thus \(|L'|\le CM\).

Where \(|Z|\ge\delta\), we claim \(|W'|\le CM\).
Self-images contribute zero to every \(\rho_j'\).
For a separated pair, \eqref{eq:selfcontained-normal} bounds
its contribution to \(W'\) by \(C(\rho_j+\rho_k)\).
A close pair containing a point away from the corner lies,
for small \(M\), in one flat wall strip.
For a vertical strip, write
\[
z=(x,y),\quad w=(u,v),\qquad
x,u\le CM,\quad y,v\ge c_\delta,
\]
and put \(e=\nabla\rho-e_x\). Then
\[
|e(z)|\le Cx^2,\qquad \|D^2\rho\|\le CM
\]
on the connecting segment.
The direct terms cancel into
\[
(e(z)-e(w))\cdot K(z-w)=O(M).
\]
The linear normal components of the vertical images cancel,
and their remaining terms are bounded by
\[
C\frac{x^2+u^2}{x+u}=O(M).
\]
The horizontal and origin images have bounded derivatives
and normal components \(O(x)\), \(O(u)\).
The horizontal-strip argument is identical.
When \(|Z|\ge\delta\), both positive vortices are away from
the corner, so every pair is covered by one of these estimates.
Therefore \(|W'|\le CM\).

Finally, on the cutoff transition,
\[
|\theta'(L-W)|
\le C|\dot Z|\,r^2/h
\le Cr^2/h^2
\le CM^{2\kappa-2}\le CM.
\]
Combining the three regions proves the lemma.
\end{proof}

\begin{remark}[(Geometric interpretation)]
The quantity \(V\) is a corrected signed moment of boundary
distance:
\[
V=(1-\theta(Z))(a\rho_1+b\rho_2-cs)
+\theta(Z)(Ah-cs),
\]
where \(\theta=1\) near the corner and \(\theta=0\) outside
a slightly larger neighborhood. Near a wall and away from
the corner, \(\rho\) is asymptotic to the normal distance
to that wall. Thus the first expression approximates the
signed normal moment of circulation. Near the corner,
the positive pair is instead represented by its effective
vortex of circulation \(A\) at \(Z\), removing the effect
of its rapid internal motion. The uniform relation
\(V=Ah-cs+o(M)\) shows that this modification preserves
the leading-order balance of the weighted distances,
even along unbounded trajectories.

The estimate \(|V'|\le CM\) expresses slow variation of
this signed moment when all boundary distances are small,
despite potentially large individual vortex velocities.
In particular, on any interval where \(M\le\varepsilon\),
\[
|V(t_2)-V(t_1)|
\le C\int_{t_1}^{t_2}M(t)\,dt
\le C\varepsilon(t_2-t_1).
\]
Because \(V\) is signed, its contributions may cancel;
the estimate therefore supplies a controlled moment
which must be combined with the effective energy to
obtain a positive boundary-exclusion barrier.
\end{remark}

\subsection{Boundary exclusion}

\begin{theorem}[(Finite-time boundary exclusion for three vortices)]
\label{thm:three-boundary-exclusion-unbounded}
Let \(\Gamma_1,\Gamma_2,\Gamma_3\ne0\), and let
\(z_1,z_2,z_3\) solve \eqref{eq:selfcontained-quadrant}
with pairwise distinct positions in \(\Omega\) for
\(0\le t<T<\infty\). Then
\[
\inf_{0\le t<T}\min_{1\le j\le3}
\operatorname{dist}(z_j(t),\partial\Omega)>0.
\]
No boundedness or trajectory-convergence assumption is required.
\end{theorem}

\begin{proof}
The same-sign case follows from
Lemma~\ref{lem:three-forced-group}.
Use the mixed-sign notation above.

We first prove
\begin{equation}
\liminf_{t\uparrow T}M(t)>0.
\label{eq:three-no-small-M}
\end{equation}

\emph{Case \(A\ne c\).}
Let \(V\) be supplied by
Lemma~\ref{lem:three-modified-moment}, and set
\[
B_0=\frac{A^2+c^2}{4\pi},\qquad
S=e^{\mathcal E/B_0}.
\]
Extend \(S=0\) at \(\ell=0\). Since \(G\ge0\), \(S\le2M\).

Choose
$
0<\beta<\frac18\min\{A,c,|A-c|\}.
$
For small \(M\), whenever \(|V|\le\beta M\), the relation
\(V=Ah-cs+o(M)\) implies
\[
h\asymp s\asymp M,\qquad |h-s|\ge c_1M.
\]
Thus \(\ell\ge|h-s|\ge c_1M\), the effective Green function
is bounded, and
\[
S\ge c_2M,\qquad
|\mathcal E'|\le Cr^2/M^4
\le CM^{2\kappa-4}\le C.
\]
In particular, \(S'\ge-CS\) in this region.

Choose \(2\lambda<\beta\) and put
$
F=\max\{|V|,\lambda S\}.
$
For all sufficiently small \(M\),
\[
c_3M\le F\le C_3M,\qquad F'\ge-CF
\quad\text{almost everywhere}.
\]
Indeed, on the \(S\)-branch,
\(|V|\le\lambda S\le2\lambda M<\beta M\), so the
effective-energy estimate applies. On the other branch,
\[
(|V|)'\ge-|V'|\ge-CM\ge-C'F.
\]
At \(\ell=0\), \(h=s=M\) and
\(V=(A-c)h+o(M)\ne0\). Thus \(F=|V|\) locally, and
\(F\) remains locally Lipschitz through effective coincidences.

Fix a small threshold \(\varepsilon>0\).
Integrate on each component of \(\{M<\varepsilon\}\).
Every entrance after \(t_0\) has
\(F\ge c_3\varepsilon\); a component beginning at \(t_0\)
has a positive initial value.
Since \(T<\infty\), Gronwall's inequality gives a uniform
positive lower bound on all such components.
Because \(F\le C_3M\), this proves
\eqref{eq:three-no-small-M}.

\emph{Case \(A=c\).}
The moments
$
P=A(X-x_3),\; Q=A(Y-y_3)
$
are respectively nondecreasing and nonincreasing.
Hence \(\ell=\sqrt{P^2+Q^2}/A\) has a limit in
\([0,\infty]\).

If this limit is positive, \(\ell\ge\delta_1>0\)
sufficiently late. Whenever \(h+s\) is small, the normal-field
and effective-pair estimates give
\[
|h'|\le Ch,\qquad |s'|\le Cs.
\]
The conditional scalar barrier for \(h+s\) proves
\eqref{eq:three-no-small-M}.

If \(\ell\to0\), monotonicity gives
$
x_3\ge X,\; Y\ge y_3.
$
Corner exclusion and
$
R\le C_*|Z|+\ell
$
imply \(|Z|\ge c_0>0\) sufficiently late.
Set \(f=\min\{Y,x_3\}\). The ordering and
\eqref{eq:selfcontained-rho} give
\[
M\le f\le\sqrt2\,(M+\ell).
\]
Thus a sequence \(M\to0\) would give \(f\to0\).

Whenever \(f\) is small, the neutral estimate applies.
If \(f=Y\), then \(X,x_3\ge c_1>0\), and
\[
V^{\mathbb H}_{Z,y}
=
\frac{A(x_3-X)}{2\pi}
\left(\ell^{-2}-|Z-\bar z_3|^{-2}\right)\ge0.
\]
The remaining images contribute \(O(Y)\), while
\[
|\mathcal R_Z|\le Cr^2/d^3\le Cd\le Ch\le CY.
\]
Hence \(Y'\ge-CY\).

If \(f=x_3\), then \(Y,y_3\ge c_1>0\), and
\[
V^{\mathrm{right}}_{3,x}
=
\frac{A(Y-y_3)}{2\pi}
\left(\ell^{-2}-|z_3+\bar Z|^{-2}\right)\ge0.
\]
The far-wall contribution is \(O(x_3)\), and
\[
|\mathcal R_{3,x}|
\le Cx_3r^2/\ell^4\le Cx_3.
\]
Hence \(x_3'\ge-Cx_3\).
The far-wall estimates in both cases follow from vanishing
normal components and uniformly bounded spatial derivatives.

Therefore \(f'\ge-Cf\) whenever \(f\) is small.
The scalar barrier excludes \(f\to0\) along a sequence,
proving \eqref{eq:three-no-small-M} also in this case.

\emph{Reduction of individual boundary approach.}
Suppose some vortex nevertheless approaches the boundary.
There must be a sequence on which an opposite-sign pair has
both \(\rho\)-values tending to zero.
Otherwise
\[
\liminf_{t\uparrow T}
\bigl(\rho_3+\min\{\rho_1,\rho_2\}\bigr)>0.
\]
Whenever \(\rho_3\) is small, both positive vortices would
then be uniformly separated from it. Equations
\eqref{eq:selfcontained-normal}--%
\eqref{eq:selfcontained-rho-derivative} give the conditional
barrier \(\rho_3'\ge-C\rho_3\).
Once \(\rho_3\) is bounded away from zero,
Lemma~\ref{lem:three-forced-group} protects the positive
pair as well, a contradiction.

Relabel the positive vortices so that
\[
\rho_1(t_n),\rho_3(t_n)\to0.
\]
If \(\rho_2(t_n)\not\to0\), pass to a subsequence with
\(\rho_2(t_n)\ge\eta>0\).
Then \eqref{eq:selfcontained-G-bound} gives
\[
G(z_1(t_n),z_2(t_n)),
\quad G(z_3(t_n),z_2(t_n))\to0.
\]
Indeed, choosing the wall nearest a small-\(\rho\) target,
whose normal coordinate is \(m\to0\), gives
\[
G(z,z_2)\le
\frac1{4\pi}\log\left(1+\frac{16m}{\eta}\right)
\]
for sufficiently large \(n\).
Since the interaction between vortices \(1\) and \(3\)
is nonpositive, conservation of energy forces
\(\rho_2(t_n)\to\infty\).
This contradicts Lemma~\ref{lem:three-cutoff}.
Thus all three \(\rho_j(t_n)\to0\), and
Lemma~\ref{lem:three-effective} gives \(M(t_n)\to0\),
contradicting \eqref{eq:three-no-small-M}.

Every \(\rho_j\) therefore has a positive lower limit at \(T\).
Continuity on initial compact time intervals, together with
\(\operatorname{dist}(z_j,\partial\Omega)\ge\rho_j\),
completes the proof.
\end{proof}


\medskip


\section{Rigorous Existence of Finite-Time Interior Collapses in the Quadrant}
\label{sec:rigorous_existence_collapse}

The equations of motion for $N=3$ point vortices at $z_j(t) = x_j(t) + i y_j(t) \in \Omega$ with circulations $\Gamma_j$ are given by:
$$\frac{d\bar{z}_j}{dt} = \frac{1}{2\pi i} \sum_{\substack{k=1 \\ k \neq j}}^3 \frac{\Gamma_k}{z_j - z_k} + R_j(z_1, z_2, z_3; \mathbf{\Gamma}),$$
where $R_j(\mathbf{z})$ represents the regular velocity field induced by the image vortices reflecting across $x=0$ and $y=0$:
\[
\begin{aligned}
R_j(\mathbf z;\boldsymbol\Gamma)
=\frac{1}{2\pi i}\Bigg[
&\Gamma_j\left(
-\frac{1}{z_j-\bar z_j}
-\frac{1}{z_j+\bar z_j}
+\frac{1}{2z_j}
\right)\\
&+\sum_{\substack{k=1\\k\ne j}}^3
\Gamma_k\left(
-\frac{1}{z_j-\bar z_k}
-\frac{1}{z_j+\bar z_k}
+\frac{1}{z_j+z_k}
\right)
\Bigg].
\end{aligned}
\]
Equivalently, including self-images in the sum,
\[
R_j(\mathbf z;\boldsymbol\Gamma)
=
\frac{1}{2\pi i}
\sum_{k=1}^3\Gamma_k
\left(
-\frac{1}{z_j-\bar z_k}
-\frac{1}{z_j+\bar z_k}
+\frac{1}{z_j+z_k}
\right).
\]
We first prove a lemma for the boundary regularity.

\begin{lemma}\label{lem:boun_reg}
Fix the circulations and \(z_*=x_*+iy_*\in\Omega\). Set
\[
r_0=\frac12\min\{x_*,y_*\},
\qquad
K=\overline{D_{r_0}(z_*)}^{\,3}\subset\Omega^3.
\]
Each \(R_j\) is real-analytic on an open neighborhood of \(K\).
Moreover, there exist \(M_0,M_1<\infty\) such that, for
\(j=1,2,3\),
\[
\sup_{\mathbf z\in K}|R_j(\mathbf z)|\le M_0,
\qquad
\sup_{\mathbf z\in K}\|DR_j(\mathbf z)\|\le M_1,
\]
where \(DR_j\) is the real Jacobian with respect to the six
position coordinates.
\end{lemma}

\begin{proof}
For \(\mathbf z\in K\), every coordinate satisfies
\[
x_j\ge x_*-r_0\ge x_*/2,
\qquad
y_j\ge y_*-r_0\ge y_*/2.
\]
Consequently, for all \(j,k\), including \(j=k\),
\[
|z_j-\bar z_k|\ge y_j+y_k\ge y_*,
\qquad
|z_j+\bar z_k|\ge x_j+x_k\ge x_*,
\]
and
\[
|z_j+z_k|\ge x_j+x_k\ge x_*.
\]
Thus every denominator in \(R_j\) is uniformly separated
from zero on \(K\). Since \(R_j\) is a rational function
of the real position coordinates with no vanishing
denominators there, it is real-analytic on an open
neighborhood of \(K\). Compactness and continuity of
\(R_j\) and \(DR_j\) give the asserted bounds, uniformly
over the three indices.
\end{proof}

Then, we prove the theorem for the existence of finite-time interior 3-vortex collapse in the quadrant.

\begin{theorem}[(Finite-time interior collapse for prescribed circulations)]
\label{thm:collapse}
Let \(\Gamma_1,\Gamma_2,\Gamma_3\in\mathbb R\setminus\{0\}\)
satisfy
\[
\Gamma_1\Gamma_2+\Gamma_1\Gamma_3+\Gamma_2\Gamma_3=0.
\]
There exist a noncollinear profile
\(\mathbf A=(A_1,A_2,A_3)\), with pairwise distinct components,
and \(\mu\in\mathbb R\setminus\{0\}\), depending only on the
circulations, with the following property.

For every prescribed collapse point
\(z_*=x_*+iy_*\in\Omega\), there exist constants \(T_0,C>0\)
such that, for every prescribed collapse time
\(t_*\in(0,T_0)\), the quadrant vortex equations admit
a smooth solution with pairwise distinct positions on
\([0,t_*)\), satisfying 
\[
\left|
z_j(t)-z_*-(t_*-t)^{1/2+i\mu}A_j
\right|
\le C(t_*-t),
\qquad j=1,2,3.
\]
The entire solution lies in the disc \(D_{r_0}(z_*)\), where
$
r_0=\frac12\min\{x_*,y_*\}.
$
In particular, the three vortices collapse simultaneously
at the prescribed interior point \(z_*\).
\end{theorem}

\begin{proof}
Define the actual complex planar velocity by
\[
P_j(\mathbf z)
=
\frac{i}{2\pi}\sum_{k\ne j}
\frac{\Gamma_k}{\bar z_j-\bar z_k}.
\]
The quadrant equations are
\begin{equation}
\dot z_j=P_j(\mathbf z)+\overline{R_j(\mathbf z)},
\label{eq:ssm-quadrant}
\end{equation}
where
\[
R_j(\mathbf z)
=
\frac1{2\pi i}\sum_{k=1}^3\Gamma_k
\left(
-\frac1{z_j-\bar z_k}
-\frac1{z_j+\bar z_k}
+\frac1{z_j+z_k}
\right).
\]
On \(\overline{D_{r_0}(z_*)}^{\,3}\), the coordinates satisfy
\(x_j\ge x_*/2\), \(y_j\ge y_*/2\). Hence
\[
|z_j-\bar z_k|\ge y_*,
\qquad
|z_j+\bar z_k|,\ |z_j+z_k|\ge x_*.
\]
Thus \(R\), regarded as a function of the six real position
coordinates, is smooth in an open neighborhood of
\((z_*,z_*,z_*)\), including configurations where the real
vortices coincide. Its derivatives are bounded on a smaller
compact neighborhood. Only the direct real-vortex interaction
is singular at this diagonal configuration.


\medskip
\noindent
\emph{1. Construction of a planar collapse profile.}

The circulation condition implies
\[
\gamma^2=\sum_j\Gamma_j^2>0,
\qquad \gamma=\sum_j\Gamma_j.
\]
After relabeling, suppose first that
\[
(\Gamma_1,\Gamma_2,\Gamma_3)=(a,b,-c),
\qquad a,b>0.
\]
Then \(c=ab/S\), where \(S=a+b\). Set
\[
\eta=\sqrt{\frac bS},
\qquad
(q_1,q_2,q_3)=(0,1,1+i\eta),
\qquad
m=\frac{b-c(1+i\eta)}{S-c},
\qquad p_j=q_j-m.
\]
The triangle is noncollinear, and
\[
S-c=\frac{a^2+ab+b^2}{S}>0,
\qquad
\sum_j\Gamma_jp_j=0.
\]
Define
\[
d=\frac{c}{2\pi\eta(1+\eta^2)}>0,
\qquad
\omega=\frac1{2\pi}
\left(S-\frac{c}{1+\eta^2}\right)>0,
\qquad \kappa=-d+i\omega.
\]
Using \(a\eta^2=c\), direct calculation gives the relations of the actual complex planar velocity $P_j$,
\[
P_2(\mathbf q)-P_1(\mathbf q)=\kappa,
\qquad
P_3(\mathbf q)-P_2(\mathbf q)=i\eta\,\kappa.
\]
Therefore \(P_j(\mathbf q)-\kappa q_j\) is independent of \(j\).
Since \(\sum_j\Gamma_jP_j(\mathbf q)=0\), this common value
is \(-\kappa m\). Translation invariance yields
\(P_j(\mathbf p)=\kappa p_j\).

Consequently, with
\[
A_j=\sqrt{2d}\,p_j,
\qquad
\mu=-\frac{\omega}{2d}\ne0,
\qquad
\beta=\frac12+i\mu,
\]
homogeneity gives
\begin{equation}
P(\mathbf A)=-\beta\mathbf A.
\label{eq:ssm-profile}
\end{equation}
If the original circulations are the negatives of
\((a,b,-c)\), conjugate this profile and replace \(\mu\)
by \(-\mu\), using
\[
P_{-\boldsymbol\Gamma}(\overline{\mathbf A})
=
\overline{P_{\boldsymbol\Gamma}(\mathbf A)}.
\]
Thus \eqref{eq:ssm-profile} holds for every circulation
triple in the theorem.

\medskip
\noindent
\emph{2. The rescaled spectrum.}

Define
\[
F(\mathbf V)=\beta\mathbf V+P(\mathbf V),
\qquad L=DF(\mathbf A),
\]
where derivatives are taken in real coordinates.
Then \(F(\mathbf A)=0\).

By Appendix A of \cite{Yang2}
\begin{equation}
\sigma(L)
=
\left\{
0,0,1,1,\frac12+i\mu,\frac12-i\mu
\right\},
\label{eq:ssm-shape-spectrum}
\end{equation}
with algebraic multiplicities.

\medskip
\noindent
\emph{3. A smooth autonomous augmentation.}

Set \(s=t_*-t\), \(\tau=-\log s\),
\(\beta=\frac12+i\mu\), and \(w=e^{-\beta\tau}\).
Then \(w'=-\beta w\), primes denote a new time
variable \(\tau\) and \(dt/d\tau=s=|w|^2\).
Writing \(\mathbf z=z_*\mathbf1+w\mathbf V\), and using
\[
P(z_*\mathbf1+w\mathbf V)=\frac{P(\mathbf V)}{\bar w},
\]
the physical equation
\(\dot{\mathbf z}=P(\mathbf z)+\overline{R(\mathbf z)}\)
gives
\[
\begin{aligned}
w(\mathbf V'-\beta\mathbf V)
&=\frac{d\mathbf z}{d\tau}\\
&=|w|^2
\left[
\frac{P(\mathbf V)}{\bar w}
+\overline{R(z_*\mathbf1+w\mathbf V)}
\right]\\
&=wP(\mathbf V)
+|w|^2\overline{R(z_*\mathbf1+w\mathbf V)}.
\end{aligned}
\]
Dividing by \(w\ne0\) and setting
\(F(\mathbf V)=\beta\mathbf V+P(\mathbf V)\), we obtain the autonomous augmented system:
\begin{equation}
\begin{aligned}
\mathbf V'
&=
F(\mathbf V)
+\bar w\,
\overline{R(z_*\mathbf1+w\mathbf V)},\\
w'&=-\beta w,
\end{aligned}
\label{eq:ssm-augmented}
\end{equation}
Since \(R\) is smooth near the interior diagonal and \(P\)
is smooth near the collision-free profile \(\mathbf A\),
the right-hand side extends smoothly to \(w=0\).

For \(\mathbf V\) near \(\mathbf A\) and \(|w|\) small,
all arguments of \(R\) lie near the interior diagonal.
The augmented vector field is therefore smooth as a real
eight-dimensional vector field, including at
$
p_*=(\mathbf A,0).
$
This point is an equilibrium.

At the interior diagonal, the image-induced velocities agree;
write \(U_*=\overline{R_j(z_*\mathbf1)}=u_*+iv_*\),
independently of \(j\). The product rule gives
\[
D\!\left[\bar w\,\overline{R(z_*\mathbf1+w\mathbf V)}\right]
_{(\mathbf A,0)}[\mathbf v,\zeta]
=\bar\zeta\,U_*\mathbf1,
\]
because all terms involving derivatives of \(R\) are
multiplied by \(\bar w=0\). Thus, using real and imaginary
coordinates for each complex component, the augmented
Jacobian is
\[
\mathcal L=
\begin{pmatrix}
L&B_*\\
0&S_*
\end{pmatrix},
\qquad
B_*=
\begin{pmatrix}M_*\\M_*\\M_*\end{pmatrix},
\qquad
M_*=
\begin{pmatrix}
u_*&v_*\\
v_*&-u_*
\end{pmatrix},
\qquad
S_*=
\begin{pmatrix}
-\frac12&\mu\\
-\mu&-\frac12
\end{pmatrix}.
\]
Here \(M_*\) represents the real-linear map
\(\zeta\mapsto U_*\bar\zeta\), while the zero lower-left
block follows from the independence of \(w'=-\beta w\)
from \(\mathbf V\). Consequently,
\[
\det(\lambda I_8-\mathcal L)
=
\det(\lambda I_6-L)
\left[\left(\lambda+\frac12\right)^2+\mu^2\right].
\]

Thus
\[
\sigma(S_*)=
\left\{-\frac12+i\mu,-\frac12-i\mu\right\},
\]
while every eigenvalue of \(L\) has nonnegative real part.

Apply the local strong stable manifold theorem, using any
spectral threshold \(-\alpha\) with \(0<\alpha<1/2\).
It gives a real two-dimensional local strong stable manifold
\(W^{ss}_{\mathrm{loc}}(p_*)\), tangent to the stable spectral
subspace of \(\mathcal L\).
The zero eigenvalues in the complementary subspace do not
obstruct this theorem.

The projection of the stable tangent space onto the
\(w\)-plane is an isomorphism. Indeed, a stable tangent
vector with \(w=0\) would give an exponentially decaying
solution of \(\mathbf v'=L\mathbf v\), which is impossible
unless \(\mathbf v=0\), by \eqref{eq:ssm-shape-spectrum}.
The dimensions are both two.

The inverse function theorem therefore represents the
stable manifold, after shrinking its neighborhood, as
\begin{equation}
\mathbf V=\mathbf A+h(w),
\qquad |w|<\varepsilon,
\label{eq:ssm-stable-graph}
\end{equation}
where \(h\) is \(C^1\) as a real map and \(h(0)=0\).
In particular,
\begin{equation}
|h(w)|\le C_h|w|.
\label{eq:ssm-graph-bound}
\end{equation}
No assertion that \(Dh(0)=0\) is needed.

Local invariance gives the graph equation
\[
Dh(w)[-\beta w]
=
F(\mathbf A+h(w))
+\bar w\,
\overline{
R(z_*\mathbf1+w(\mathbf A+h(w)))
}.
\]
Consequently, for any sufficiently small \(w_0\ne0\),
\[
w(\tau)=w_0e^{-\beta\tau},
\qquad
\mathbf V(\tau)=\mathbf A+h(w(\tau))
\]
solves \eqref{eq:ssm-augmented} for every \(\tau\ge0\):
the scale remains in the graph domain because
\(|w(\tau)|=|w_0|e^{-\tau/2}\).

\medskip
\noindent
\emph{4. Reconstruction of the collapsing solution.}

Choose \(\varepsilon>0\) so that the stable graph is defined
for \(|w|<\varepsilon\), with
\[
|h(w)|\le C_h|w|,
\qquad C_h\ge1.
\]
Set
\[
A_{\max}=\max_j|A_j|,
\qquad
d_{\min}=\min_{j<k}|A_j-A_k|>0,
\]
and define
\[
T_0=
\min\left\{
\varepsilon^2,\,
\frac{r_0^2}{4A_{\max}^2},\,
\frac{r_0}{2C_h},\,
\frac{d_{\min}^2}{16C_h^2}
\right\}>0.
\]
For every \(0<t_*<T_0\), the scale satisfies
\[
|w(\tau)|=\sqrt{t_*}\,e^{-\tau/2}<\varepsilon
\qquad(\tau\ge0),
\]
so the trajectory remains in the stable graph domain.
Moreover,
\[
|z_j(t)-z_*|
\le A_{\max}\sqrt{t_*}+C_ht_*<r_0,
\]
and
\[
|V_j(\tau)-V_k(\tau)|
\ge d_{\min}-2C_h\sqrt{t_*}
>\frac{d_{\min}}2>0.
\]
Thus this single choice of \(T_0\), independent of \(t_*\),
ensures both interior containment and pairwise distinctness
throughout every constructed solution.

For every \(0<t_*<T_0\), set
\[
w_0=t_*^\beta,\qquad
t=t_*(1-e^{-\tau}),\qquad
\mathbf z=z_*\mathbf1+w(\tau)\mathbf V(\tau).
\]
Then
\[
w(\tau)=(t_*-t)^\beta,
\qquad
\frac{dt}{d\tau}=t_*e^{-\tau}=|w(\tau)|^2.
\]
Using \eqref{eq:ssm-augmented},
\[
\frac{d\mathbf z}{d\tau}
=
wP(\mathbf V)+|w|^2\overline{R(\mathbf z)}.
\]
Translation invariance and complex scaling of the planar
field give
\[
P(z_*\mathbf1+w\mathbf V)=\frac{P(\mathbf V)}{\bar w}.
\]
Dividing by \(dt/d\tau=|w|^2\) therefore yields exactly
\eqref{eq:ssm-quadrant}.

By \eqref{eq:ssm-graph-bound},
\[
\begin{aligned}
z_j(t)
&=z_*+wA_j+w\,h_j(w)\\
&=z_*+(t_*-t)^{1/2+i\mu}A_j+O(t_*-t).
\end{aligned}
\]
The bound is uniform for \(0<t_*<T_0\).

Finally,
\[
|z_j(t)-z_*|
\le
\sqrt{t_*}\max_k|A_k|+C_ht_*<r_0
\]
when \(0<t_*<T_0\).
Likewise,
\[
|V_j(\tau)-V_k(\tau)|
\ge |A_j-A_k|-2C_h\sqrt{t_*}>0
\qquad(j\ne k).
\]
Thus the physical vortices remain distinct for \(t<t_*\)
and stay strictly inside the quadrant.
Their trajectories are smooth because they solve the smooth
physical vector field away from collisions.
The asymptotic expansion shows that all three converge to
\(z_*\) as \(t\uparrow t_*\).
\end{proof}

\begin{remark}[(Dimension of the strong stable manifold)]
The two-dimensional strong stable manifold arises from using
the complex scale \(w\), whose equation has the stable
eigenvalues \(-\frac12\pm i\mu\). Writing \(w=re^{i\theta}\)
gives
\[
r'=-\frac12r,\qquad \theta'=-\mu,
\qquad
\Psi=\theta-2\mu\log r=\text{constant}\pmod{2\pi}.
\]
Thus the manifold includes both motion along a contracting
orbit and variation of its logarithmic phase offset \(\Psi\).
In the real-scale half-plane formulation \cite{Yang2}, this phase is fixed
through the chosen profile and reconstruction, and only
\(\eta'=-\eta/2\) is adjoined. That augmented system therefore
has a one-dimensional strong stable manifold.

Our choice \(w(0)=t_*^{1/2+i\mu}\) fixes \(\Psi=0\), selecting
one phase branch. However, a fixed-phase branch is the
logarithmic spiral
$
w=r\,e^{i(2\mu\log r+\Psi)},
$
which does not extend to a regular \(C^1\) curve through
\(w=0\) when \(\mu\ne0\). Consequently, the two augmentations
are not related by a regular \(C^1\) coordinate change at
collapse. Their different stable-manifold dimensions reflect
different treatments of the phase, while both constructions
produce the required collapsing solutions.
\end{remark}

\begin{remark}[(Asymptotic self-similarity)]
The collapsing solutions are asymptotically self-similar.
Indeed, writing \(s=t_*-t\), the theorem gives
\[
\frac{e^{-i\mu\log s}}{\sqrt{s}}
\bigl(z_j(t)-z_*\bigr)
=
A_j+O(\sqrt{s})
\longrightarrow A_j.
\]
Thus, after translating the collapse point to the origin,
rescaling lengths by \(s^{-1/2}\), and removing the
logarithmic rotation, the vortex configuration converges
to the fixed nondegenerate triangle \(\mathbf A\).

In particular, for \(j\ne k\),
\[
|z_j(t)-z_k(t)|
=
\sqrt{s}\left(|A_j-A_k|+O(\sqrt{s})\right).
\]
Consequently, all side-length ratios and interior angles
converge to those of the profile triangle. The characteristic
length scale is \(\sqrt{t_*-t}\), while the leading rotation
angle is \(\mu\log(t_*-t)\). Since \(\mu\ne0\), the triangle
undergoes infinitely many rotations as it collapses.
\end{remark}


\medskip


\end{document}